\input amstex
\documentstyle{amsppt}
\topmatter
\nologo
\magnification=\magstep1
\pagewidth{5.1 in}
\pageheight{6.7 in}
\abovedisplayskip=10pt
\belowdisplayskip=10pt
\title ON RAMANUJAN'S CUBIC CONTINUED FRACTION AND
EXPLICIT EVALUATIONS OF THETA-FUNCTIONS
\endtitle

\author Chandrashekra ADIGA, Taekyun KIM,$^*$ \\
M. S. Mahadeva NAIKA, and H. S. MADHUSUDHAN \endauthor

\affil{\it
Department of Studies in Mathematics, University of Mysore \\
Manasagangotri, MYSORE-570 006, INDIA
\\
$^*$Institute of Science Education, Kongju National University \\
Kongju 314-701, REPUBLIC OF KOREA
}\endaffil

\abstract{In this paper we give two integral representations for the Ramanujan's cubic continued fraction $V(q)$ and
also derive a modular equation relating $V(q)$ and $V(q^3).$  We also establish some modular equations and a transformation
formula for Ramanujan's theta-function.
As an application of these, we compute several new explicit evaluations of theta-functions and Ramanujan's cubic continued
fraction.}\endabstract
\footnote""{Key words. Cubic continued fraction, integral representation, modular equation, theta-function,
transformation formula.}
\footnote""{2000 AMS Subject Classification: 11A55, 33D15, 33D20.}
\leftheadtext{}
\rightheadtext{}
\endtopmatter
\document
\TagsOnRight
\noindent

\head\bf 1. Introduction \endhead

The celebrated Rogers-Ramanujan continued fraction is defined by
$$R(q):=\frac{q^{1/5}}1 \frac{q}{+~1}\frac{q^2}{+~1}\frac{q^3}{+~1~+\cdots}, \quad |q|<1.\tag1.1$$
On page 46 in his `lost' notebook [19, p. 46], Ramanujan claims that
$$\align
R(q)&=\frac{\sqrt5-1}{2}\exp\left((-1/5)\int_q^1\frac{(1-t)^5(1-t^2)^5\cdots}{(1-t^5)(1-t^{10})\cdots}\frac{dt}{t}\right)
\tag1.2 \\
&=\frac{\sqrt5-1}{2}-\frac{\sqrt5}{1+\frac{3+\sqrt5}{2}}\exp\left((1/5)\int_0^q
\frac{(1-t)^5(1-t^2)^5\cdots}{(1-t^{1/5})(1-t^{2/5})\cdots}\frac{dt}{t^{4/5}}\right),
\tag1.3\endalign$$
where $0 < q < 1.$ The first equality (1.2) was proved by G.E. Andrews [3] and the second equality (1.3)
was proved by Seung Hwan Son [20]. On page 207 of his 'lost' notebook Ramanujan also recorded six identities
involving integrals of theta-functions.  All these identities were proved by S.H. Son [20] and were generalized
by Chandrashekar Adiga, K. R. Vasuki and M. S. Mahadeva Naika [1]. On page 365 of his 'lost' notebook [19], Ramanujan
wrote five modular equations relating $R(q)$ with $R(-q),R(q^2),R(q^3),R(q^4),$
and $R(q^5).$

Ramanujan eventually found several generalizations and ramifications of (1.1) which are recorded in his 'lost' notebook.
These and related works may be found in the papers by S. Bhargava [8], S. Bhargava
and C. Adiga [9], [10], R. Y. Denis [13], [14], [15].

On page 366 of his 'lost' notebook [19], Ramanujan investigated the continued fraction
$$V(q):=\frac{q^{1/3}}1 \frac{q+q^2}{+~1}\frac{q^2+q^4}{+~1~+\cdots}, \quad |q|<1,\tag1.4$$
which is known as Ramanujan's cubic continued fraction.  Analogous to those relations for $R(q)$ which are mentioned
above, H. H. Chan [12] established several modular equations relating $V(q)$ with $V(-q),V(q^2)$ and $V(q^3).$
An example of these modular equations is
$$V^3(q)=V(q^3)\frac{1-V(q^3)+V^2(q^3)}{1+2V(q^3)+4V^2(q^3)}.\tag1.5$$
In Section 2, we will establish two integral representations for $V(q),$ which are analogous to (1.2) and (1.3).
We also give a simple proof of (1.5).

In Ramanujan's theory of theta-functions, the three theta-functions that play central roles are defined by
$$\varphi(q):=\sum_{n=-\infty}^{\infty}q^{n^2}=\frac{(-q;q)_\infty}{(q;-q)_\infty},\tag1.6$$
$$\psi(q):=\sum_{n=0}^{\infty}q^{n(n+1)/2}=\frac{(q^2;q^2)_\infty}{(q;q^2)_\infty},\tag1.7$$
and
$$f(-q):=\sum_{n=-\infty}^{\infty}(-1)^{n}q^{n(3n-1)/2}=(q;q)_\infty,\tag1.8$$
where
$$(a;q)_\infty:=\prod_{n=0}^\infty(1-aq^n),\quad|q|<1.$$

In Chapter 16 of his second notebook [18], Ramanujan records many transformation formulas for $\varphi(q), \psi(q)$
and $f(-q).$ Four of the most important
transformation formulas are given by
$$\sqrt\alpha\varphi(e^{-\alpha^2})=\sqrt\beta\varphi(e^{-\beta^2}),\quad\alpha\beta=\pi,
\tag1.9$$
$$2\sqrt\alpha\psi(e^{-2\alpha^2})=\sqrt\beta e^{\alpha^2/4}\varphi(-e^{-\beta^2}),\quad
\alpha\beta=\pi,\tag1.10$$
$$e^{-\alpha/12}\root4\of{\alpha} f(-e^{-2\alpha})=e^{-\beta/12}\root4\of{\beta}f(e^{-2\beta}),
\quad\alpha\beta=\pi^2,\tag1.11$$
and
$$e^{-\alpha/24}\root4\of{\alpha} f(e^{-\alpha})=e^{-\beta/24}\root4\of{\beta}f(e^{-\beta}),
\quad\alpha\beta=\pi^2.\tag1.12$$

In Section 3, we will prove the transformation formula
$$\root4\of{\alpha}
e^{-\alpha/8}\psi(-e^{-\alpha})=e^{-\beta/8}\root4\of{\beta}
\psi(-e^{-\beta}). \tag1.13$$
On page 204 of his second notebook [18, p. 204], Ramanujan claims that
$$\left(\frac{\sqrt5+1}2+R(e^{-2\pi\alpha})\right)
 \left(\frac{\sqrt5+1}2+R(e^{-2\pi\beta})\right)=\frac{5+\sqrt5}2\tag1.14$$
and
$$\left(\frac{\sqrt5-1}2-R(-e^{-2\pi\alpha})\right)
 \left(\frac{\sqrt5-1}2-R(-e^{-2\pi\beta})\right)=\frac{5-\sqrt5}2\tag1.15$$
Identities (1.14) and (1.15) were first proved by G. N. Watson [21].  H. H. Chan [12] has
proved several identities for $V(q)$ which are similar to (1.14) and (1.15).
In Section 4, we will establish three reciprocity
theorems for $V(q)$ which are also similar to (1.14) and (1.15).

Ramunajan has recorded many modular equations in his notebooks [5, pp. 204-237], [6, pp. 156-160] which are very
useful in the computation of class invariants and the values of theta-functions.
In the literature not much attention has been given to find the values of $\psi(q)$ and $\varphi(q).$
But Ramanujan recorded several values of $\varphi(q)$ in his notebooks. For example
$$\varphi(e^{-\pi})=\frac{\pi^{1/4}}{\Gamma(3/4)},$$
$$\psi(e^{-\pi})=2^{-5/8}e^{\pi/8}\frac{\pi^{1/4}}{\Gamma(3/4)},$$
$$\frac{\varphi(e^{-\pi})}{\phi(e^{-3\pi})}=\root 4 \of{6\sqrt3-9},$$
and
$$\frac{\varphi(e^{-\pi})}{\varphi(e^{-5\pi})}=\sqrt{5\sqrt5-10}.$$
J. M. Borwein and P. B. Borwein [11] first observed that class invariants could be used to calculate certain
values of $\varphi(e^{-n\pi}).$
Bruce C. Berndt [6] has verified all values for $\varphi(e^{-n\pi})$ claimed by
Ramanujan by combining Ramanujan's class invariants with modular equations.

In Section 5, we derive some modular equations and briefly discuss evaluations of theta-functions $\psi(q)$
and $\varphi(q).$
At the end of this section we compute some interesting new explicit evaluations of $V(q).$
Our work is sequel to the works of B. C. Berndt et. al. [7], H. H. Chan [12], K.
Ramachandra [16], K. G. Ramanathan [17] and C. Adiga, K. R. Vasuki and M. S. Mahadeva Naika [2].
In this paper we adopt existing methods in the literature and work with $\psi(q)$ instead of $\varphi(q)$
as is done in [12].

\head \bf 2. Integral representations for $V(q)$ and modular equation relating $V(q)$ with
$V(q^3)$ \endhead

\proclaim{Theorem 2.1} We have
$$\align
V(q)&=\frac1{\root 3\of{-1+9\exp\left(\int_q^1\varphi^2(-t)\varphi^2(-t^3)
\frac{dt}{t}\right)} }\tag2.1 \\
&=\frac12 \root 3\of{1-\exp\left(-8\int_0^q\psi^2(t)\psi^2(t^3)
dt\right)}.\tag2.2
\endalign$$
\endproclaim
\demo{Proof of (2.1)}  Let $F(q):=\frac{\psi^4(q)}{q\psi^4(q^3)}.$
Using (1.7) and then taking logarithm on both sides, we find that
$$\align\log F(q)&=4\sum_{n=1}^\infty\{\log(1-q^{2n})-\log(1-q^{2n-1})\} \\
&+4\sum_{n=1}^\infty\{\log(1-q^{3(2n-1)})-\log(1-q^{6n})\}-\log q.
\endalign$$
Taking derivative on both sides of the above identity, we deduce that
$$\frac d{dq}[\log F(q)]=-\frac1q\left[1+4\sum_{n=1}^\infty\frac{(-1)^nnq^n}{1-q^n}
-12\sum_{n=1}^\infty\frac{(-1)^nnq^{3n}}{1-q^{3n}}\right].$$
Using Entry 3 (iv) of Chapter 19 of Ramanujan's notebooks [4, p. 226], we find that
$$\frac d{dq}[\log F(q)]=-\frac{\varphi^2(-q)\varphi^2(-q^3)}{q}$$
Integrating the above identity over $[q,1]$ on both sides, we obtain
$$\log F(q)=\int_q^1\varphi^2(-t)\varphi^2(-t^3)\frac{dt}t+\log9.\tag2.3$$
Using Entry 1 (i) of Chapter 20 of Ramanujan's notebooks [4, p. 345] in (2.3) and then exponentiating, we obtain (2.1).
\enddemo
\demo{Proof of (2.2)} Let $H(q):=\frac{\varphi^4(-q)}{\varphi^4(-q^3)}.$
Using (1.6) and then taking logarithm on both sides, we find that
$$\align\log H(q)&=4\sum_{n=1}^\infty\{\log(1-q^{n})-\log(1+q^{n})\} \\
&+4\sum_{n=1}^\infty\{\log(1+q^{3n})-\log(1-q^{3n})\}.
\endalign$$
Taking derivative on both sides of the above identity, we see that
$$\frac d{dq}[\log H(q)]=-\frac8q\left[\sum_{n=1}^\infty\frac{nq^n}{1-q^{2n}}
-3\sum_{n=1}^\infty\frac{nq^{3n}}{1-q^{6n}}\right].$$
Using Entry 1 (iii) of Chapter 19 of Ramanujan's notebooks [4, p. 225], we deduce that
$$\frac d{dq}[\log H(q)]=-8\psi^2(q)\psi^2(q^3).$$
Integrating the above identity over $[0,q]$ on both sides and then exponentiating, we obtain
$$H(q)=\exp\left(-8\int_0^q\psi^2(t)\psi^2(t^3)dt\right).\tag2.4$$
We deduce (2.2) on employing the following identity [4, p. 345] in (2.4):
$$1-8V^3(q)=\frac{\psi^4(-q)}{\psi^4(-q^3)}.\tag2.5$$
\enddemo
\proclaim{Theorem 2.2} Let $V(q)$ be defined as in (1.4).  Then
$$V^3(q)=V^3(q)\frac{1-V(q^3)+V^2(q^3)}{1+2V(q^3)+4V^2(q^3)}.\tag2.6$$
\endproclaim
\demo{Proof} To prove this theorem we require the following identities [4, p. 345, Entry 1 (i), (ii)]:
$$1+\frac1{V(q^3)}=\frac{\psi(q)}{q\psi(q^9)},\tag2.7$$
$$1+\frac1{V^3(q)}=\frac{\psi^4(q)}{q\psi^4(q^3)},\tag2.8$$
and
$$1-3q\frac{\psi(q^9)}{\psi(q)}=\left[1-9q\frac{\psi^4(q^3)}{\psi^4(q)}\right]^{1/3}.\tag2.9$$
Using (2.7) and (2.8) in (2.9), we obtain
$$1-\frac{3V(q^3)}{1+V(q^3)}=\left[1-\frac{9V^3(q)}{1+V^3(q)}\right]^{1/3}.\tag2.10$$
Cubing both sides of the above identity (2.10), we deduce (2.6).
\enddemo

\head \bf 3. Transformation formula for $\psi$ \endhead

In the following theorem we prove a transformation formula for $\psi$ akin to Ramanujan's transformation formulas (1.9)-(1.12).

\proclaim{Theorem 3.1} If $\alpha\beta=\pi^2,$ then
$$\root4\of{\alpha}e^{-\alpha/8}\psi(-e^{-\alpha})=e^{-\beta/8}\root4\of{\beta}\psi(-e^{-\beta}).\tag3.1$$
\endproclaim
\demo{Proof} Interchanging $\alpha$ and $\beta$ in (1.10), we find that
$$2\sqrt\beta\psi(e^{-2\beta^2})=\sqrt\alpha e^{\beta^2/4}\varphi(-e^{-\alpha^2}).\tag3.2$$
Using (1.10) and (3.2), we deduce that
$$\frac{\psi(e^{-2\alpha^2})}{\psi(e^{-2\beta^2})}=\frac{\beta}{\alpha}\frac{e^{\alpha^2/4}}{e^{\beta^2/4}}
\frac{\varphi(-e^{-\beta^2})}{\varphi(-e^{-\alpha^2})},\quad\alpha\beta=\pi.\tag3.3$$
Replacing $\alpha$ by $\sqrt\alpha$ and $\beta$ by $\sqrt\beta$ in (3.3), we find that
$$\frac{\psi(e^{-2\alpha})\psi(-e^{-\alpha})}{\psi(e^{-2\beta})\psi(-e^{-\beta})}=\sqrt{\frac\beta\alpha}
\frac{e^{\alpha/4}}{e^{\beta/4}},\quad\alpha\beta=\pi^2.\tag3.4$$
Using Entry 25 (iv) of Chapter 16 of Ramanujan's notebooks [4, p. 40] in (3.4), we obtain (3.1).
\enddemo
\noindent
{\bf Remark.} The transformation formula (3.1) can also be proved by using the identities (1.11) and (1.12).

\head\bf 4. Reciprocity theorems for $V(q)$\endhead

\proclaim{Theorem 4.1} We have
\roster
\item"(i)" If $\alpha\beta=1,$ then
$$\left[1+\frac1{V(-e^{-\pi\alpha})}\right]\left[1+\frac1{V(-e^{-\pi\beta})}\right]=3.\tag4.1$$
\item"(ii)" If $3\alpha\beta=1,$ then
$$\left[1+\frac1{V^3(-e^{-\pi\alpha})}\right]\left[1+\frac1{V^3(-e^{-\pi\beta})}\right]=9.\tag4.2$$
\item"(iii)" If $3\alpha\beta=1,$ then
$$\left[1+\frac1{V^3(-e^{-\sqrt2\pi\alpha})}\right]\left[1-8V^3(e^{-\sqrt2\pi\beta})\right]=9.\tag4.3$$
\endroster
\endproclaim
\demo{Proof of (4.1)} Using (2.7), we find that
$$\left[1+\frac1{V(-e^{-\pi\alpha})}\right]\left[1+\frac1{V(-e^{-\pi\beta})}\right]=
\frac{\psi(-e^{-\pi\alpha/3})\psi(e^{-\pi\beta/3})}{e^{(-\pi/3)(\alpha+\beta)}\psi(-e^{-3\pi\alpha})
\psi(-e^{-3\pi\beta})}.\tag4.4$$
From the transformation formula (3.1), we deduce that
$$\frac{\psi(-e^{-\pi\alpha/3})}{e^{(\pi/8)((\alpha/3)-3\beta)}\psi(-e^{-3\pi\beta})}=
\root4\of{\frac{9\beta}{\alpha}}.\tag4.5$$
Interchanging $\alpha$ and $\beta$ in (4.5), we obtain
$$\frac{\psi(-e^{-\pi\beta/3})}{e^{(\pi/8)((\beta/3)-3\alpha)}\psi(-e^{-3\pi\alpha})}=
\root4\of{\frac{9\alpha}{\beta}}.\tag4.6$$
Using (4.5) and (4.6) in (4.4), we obtain (4.1).
\enddemo
\redefine\a{\alpha}
\redefine\b{\beta}
\demo{Proof of (4.2)} Using (2.8), we find that
$$\left[1+\frac1{V^3(-e^{-\pi\alpha})}\right]\left[1+\frac1{V^3(-e^{-\pi\beta})}\right]=
\frac{\psi^4(-e^{-\pi\alpha})\psi^4(e^{-\pi\beta})}{e^{-\pi(\alpha+\beta)}\psi^4(-e^{-3\pi\alpha})
\psi^4(-e^{-3\pi\beta})}.\tag4.7$$
From the transformation formula (3.1), we deduce that
$$\frac{\psi^4(-e^{-\pi\alpha})}{e^{(\pi/2)(\alpha-3\beta)}\psi^4(-e^{-3\pi\beta})}=
\frac{3\beta}{\alpha}.\tag4.8$$
Interchanging $\a$ and $\b$ in (4.8), we obtain
$$\frac{\psi^4(-e^{-\pi\beta})}{e^{(\pi/2)(\beta-3\alpha)}\psi^4(-e^{-3\pi\alpha})}=
\frac{3\a}{\b}.\tag4.9$$
Using (4.8) and (4.9) in (4.7), we obtain (4.2).
\enddemo
\demo{Proof of (4.3)} Using (2.8), we find that
$$\left[1+\frac1{V^3(-e^{-\sqrt2\pi\alpha})}\right]=\frac{\psi^4(e^{-\sqrt2\pi\a})}
{e^{-\sqrt2\pi\a}\psi^4(e^{-3\sqrt2\pi\a})}\tag4.10$$
Using (1.10) in (4.10), we deduce that
$$\left[1+\frac1{V^3(e^{-\sqrt2\pi\alpha})}\right]=9\frac{\psi^4(-e^{-\sqrt2\pi/\a})}
{\psi^4(-e^{-\sqrt2\pi/3\a})}\tag4.11$$
Using (2.5) in (4.11), we find that
$$\left[1+\frac1{V^3(-e^{-\sqrt2\pi\alpha})}\right]=\frac9{1-8V^3(e^{-\sqrt2\pi/3\a})}
=\frac9{1-8V^3(e^{-\sqrt2\pi\b})}.$$
Hence, we complete the proof.
\enddemo

\head\bf 5. Modular equations and evaluations of theta-function and
$V(q)$\endhead

\proclaim{Theorem 5.1} Let
$$P=\frac{\psi(-q)}{q^{1/4}\psi(-q^3)}\quad\text{and}\quad Q=\frac{\varphi(q)}{\varphi(q^3)}.$$
Then
$$Q^4+P^4Q^4=9+P^4.\tag5.1$$
\endproclaim
\demo{Proof} Using Entry 1 (i) of Chapter 20 of Ramanujan's notebooks [4, p. 345] and (2.5), we find that
$$\frac{\psi^4(q)}{q\psi^4(q^3)}+\frac{\varphi^4(-q)}{\varphi^4(-q^3)}=9+\frac{\psi^4(q)\varphi^4(-q)}
{q\psi^4(q^3)\varphi^4(-q^3)}.$$
Replacing $q$ by $-q$ in the above identity, we obtain (5.1).
\enddemo
\proclaim{Theorem 5.2} Let
$$P=\frac{\psi(-q)}{q\psi(-q^9)}\quad\text{and}\quad Q=\frac{\varphi(q)}{\varphi(q^9)}.$$
Then
$$Q+PQ=3+P\tag5.2$$
\endproclaim
\demo{Proof} Using Entry 1 (i) and (ii) of Chapter 20 of Ramanujan's notebooks [4, p. 345], we deduce that
$$\frac{\psi(q^{1/3})}{q^{1/3}\psi(q^3)}+\frac{\varphi(-q^{1/3})}{\varphi(-q^3)}
=3+\frac{\psi(q^{1/3})\varphi(-q^{1/3})}{q^{1/3}\psi(q^3)\varphi(-q^3)}.$$
Replacing $q$ by $-q^3$ in the above identity, we obtain (5.2).
\enddemo
\proclaim{Theorem 5.3} Let
$$P=\frac{\psi(-q)}{q^{1/2}\psi(-q^5)}\quad\text{and}\quad Q=\frac{\varphi(q)}{\varphi(q^5)}.$$
Then
$$Q^2 +P^2Q^2=5+P^2.\tag5.3$$
\endproclaim
\demo{Proof} Changing $q$ to $-q$ in Entry 9 (iii) of Chapter 19 of Ramanujan's notebooks [4, p. 258],
we find that
$$\frac{\varphi^2(-q)}{\varphi^2(-q^5)}=1-4W,\tag5.4$$                                                                                   (5.4)
where
$$W=\frac{q\chi(-q)f(q^5)f(-q^{20})}{\varphi^2(-q^5)}.$$
Using Entry 9 (vii) of Chapter 19 of Ramanujan's notebooks [4, p. 258], Entry 10 (v) of Chapter 19 of
Ramanujan's notebooks [4, p. 262] can be written as
$$\frac{\psi^2(q)}{q\psi^2(q^5)}=1+\frac{\varphi(-q^5)f(-q^5)}{q\chi(-q)\psi^2(q^5)}.\tag5.5$$
The above identity can be written as
$$\frac{\psi^2(q)}{q\psi^2(q^5)}=1+\frac1W.\tag5.6$$
Using (5.4) and (5.6), we find that
$$\frac{\psi^2(q)}{q\psi^2(q^5)}+\frac{\varphi^2(-q)}{\varphi^2(-q^5)}=5+\frac{\psi^2(q)\varphi^2(-q)}
{q\psi^2(q^5)\varphi^2(-q^5)}.$$
Changing $q$ to $-q$ in the above identity, we obtain (5.3).
\enddemo
\proclaim{Theorem 5.4} We have
\roster
\item"(i)" $$\frac{\psi(-e^{-\pi/\sqrt5})}{e^{-\pi/2\sqrt5}\psi(-e^{-\sqrt5\pi})}=5^{1/4},\tag5.7$$
\item"(ii)" $$\frac{\psi(-e^{-3\pi/\sqrt5})}{e^{-3\pi/2\sqrt5}\psi(-e^{-3\sqrt5\pi})}
=\frac{5^{1/4}(\sqrt3+1)(\sqrt5+\sqrt3)}2,\tag5.8$$
\item"(iii)" $$\frac{\psi(-e^{-\pi/3\sqrt5})}{e^{-\pi/3\sqrt5}\psi(-e^{-3\pi/\sqrt5})}
=\frac{\sqrt3(\sqrt3-1)(\sqrt5-\sqrt3)}2,\tag5.9$$
\item"(iv)" $$\frac{\psi(-e^{-\pi/3\sqrt3})}{e^{-\pi/12\sqrt3}\psi(-e^{-\pi/\sqrt3})}
=3^{-1/4}(\root3\of4-1),\tag5.10$$
\item"(v)" $$\frac{\psi(-e^{-\pi})}{e^{-\pi/4}\psi(-e^{-3\pi})}
=\root4\of{3\sqrt3(2+\sqrt3)},\tag5.11$$
\item"(vi)" $$\frac{\psi(-e^{-\pi/3\sqrt5})}{e^{-\pi/6\sqrt5}\psi(-e^{-\sqrt5\pi/3})}
=\frac{5^{1/4}(\sqrt3-1)(\sqrt5-\sqrt3)}2,\tag5.12$$
\item"(vii)" $$\frac{\psi(-e^{-\sqrt5\pi/3})}{e^{-\sqrt5\pi/3}\psi(-e^{-3\sqrt5\pi})}
=\frac{\sqrt3(\sqrt3+1)(\sqrt5+\sqrt3)}2,\tag5.13$$
\item"(viii)" $$\frac{\psi(-e^{-\pi/\sqrt3})}{e^{-\pi/4\sqrt3}\psi(-e^{-\sqrt3\pi})}=3^{1/4},\tag5.14$$
\item"(ix)" $$\frac{\psi(-e^{-\pi\sqrt3})}{e^{-\pi\sqrt3/4}\psi(-e^{-3\sqrt3\pi})}
=\frac{\root4\of27}{\root3\of4-1},\tag5.15$$
\item"(x)" $$\frac{\psi(-e^{-\pi/3})}{e^{-\pi/3}\psi(-e^{-3\pi})}
=\sqrt3,\tag5.16$$
\item"(xi)" $$\frac{\psi(-e^{-\pi/3})}{e^{-\pi/12}\psi(-e^{-\pi})}
=\root4\of{\sqrt3(2-\sqrt3)},\tag5.17$$
\item"(xii)" $$\frac{\psi(-e^{-\pi/\sqrt3})}{e^{-\pi/\sqrt3}\psi(-e^{-3\sqrt3\pi})}
=\frac3{\root3\of{4}-1},\tag5.18$$
\item"(xiii)" $$\frac{\psi(-e^{-\pi/3\sqrt3})}{e^{-\pi/3\sqrt3}\psi(-e^{-\sqrt3\pi})}
={\root3\of{4}-1},\tag5.19$$
\endroster
\endproclaim
Proofs of the identities (5.7)-(5.19) being similar, for brevity we prove only (5.7)-(5.11).
\demo{Proof of (5.7)} Putting $\a=\frac\pi{\sqrt5}$ and $\b=\pi{\sqrt5}$ in the transformation formula
(3.1), we obtain (5.7).
\enddemo
\demo{Proof of (5.8)} By Entry 66 of Chapter 25 of Ramanujan's notebooks [5, p. 233] with $q$ replaced by $-q,$
we have
$$P^3Q^3+5PQ=Q^4-3PQ^3-3P^3Q-P^4\tag5.20$$
where
$$P=\frac{\psi(-q)}{q^{1/2}\psi(-q^5)}\quad\text{and}\quad
Q=\frac{\psi(-q^3)}{q^{3/2}\psi(-q^{15})}.$$
Using (5.7) in (5.20) with $q=e^{-\pi/\sqrt5},$ we find that
$$Q^4-5^{1/4}(3+\sqrt5)Q^3-5^{3/4}(3+\sqrt5)Q-5=0.\tag5.21$$
Putting $Q=i5^{1/4}T$ in (5.21), we deduce that
$$T^4+(3+\sqrt5)iT^3-(3+\sqrt5)iT-1=0.\tag5.22$$
The equation (5.22) can be written as
$$(T^2-1)(T^2+(3+\sqrt5)iT+1)=0.$$
Solving the above equation, we deduce that
$$T=\pm1,\quad T=\frac{-(3+\sqrt5)i\pm i\sqrt{(3+\sqrt5)^2+4}}2$$
and hence
$$Q=\pm i5^{1/4},\quad Q=\frac{5^{1/4}\left((3+\sqrt5)\pm \sqrt{(3+\sqrt5)^2+4}\right)}2.$$
Since $Q>0,$ we obtain (5.8).
\enddemo
\demo{Proof of (5.9)} Putting $q=e^{-\pi/\sqrt5}$ in (5.20) and using (5.7), we find that
$$P=5^{1/4}.\tag5.23$$
Also, using (3.1), we see that
$$Q=\frac{\root4\of{45}}{C}\tag5.24$$
where
$$C=\frac{\psi(-e^{-\pi/3\sqrt5})}{e^{-\pi/3\sqrt5}\psi(-e^{-3\pi/\sqrt5})}.$$
Using (5.23) and (5.24), we find that
$$PQ=\frac{\sqrt{15}}{C}\tag5.25$$
and
$$\frac PQ =\frac{C}{\sqrt3}.\tag5.26$$
Using (5.25) and (5.26) in (5.20), we deduce that
$$\sqrt5\left[\frac{\sqrt3}C+\frac C{\sqrt3}\right]=\left[\frac{\sqrt3}C+\frac C{\sqrt3}\right]
\left[\frac{\sqrt3}C-\frac C{\sqrt3}\right]-3\left[\frac{\sqrt3}C+\frac C{\sqrt3}\right].$$
Since $\frac{\sqrt3}C+\frac C{\sqrt3}\neq0,$ we obtain
$$\frac{\sqrt3}C-\frac C{\sqrt3}=\sqrt5+3.$$
Since $C > 0,$ solving the above equation we obtain the required result.
\enddemo
\demo{Proof of (5.10)} Let
$$P=\frac{\psi(-q)}{q^{1/4}\psi(-q^3)}\quad\text{and}\quad
Q=\frac{\psi(-q^3)}{q^{3/4}\psi(-q^{3})}.\tag5.27$$
Using (5.27) in Entry 1 (ii) of Chapter 20 of Ramanujan's notebooks [4, p. 345], we find that
$$Q^3-P^3Q^2-3P^2Q-3P=0.\tag5.28$$
Putting $q=e^{-\pi/3\sqrt3}$ in (5.27), and using the transformation formula (3.1), we see that
$$Q=3^{1/4}.\tag5.29$$
Using (5.29) in (5.28), we deduce that
$$P^3+3^{3/4}P^2+\sqrt3 P-3^{1/4}=0.\tag5.30$$
Putting $P=3^{-1/4}T$ in (5.30), we find that
$$T^3+3T^2+3T-3=0\tag5.31$$
Putting $x=T+1$ in (5.31), we obtain
$$x=\root3\of4.\tag5.32$$
Using (5.32), we obtain the required result.
\enddemo
\demo{Proof of (5.11)} Putting $q=e^{-\pi/3}$ in (5.28), we find that
$$P=\frac{\psi(-e^{-\pi/3})}{e^{-\pi/12}\psi(-e^{-\pi})}\tag5.33$$
and
$$Q=\frac{\psi(-e^{-\pi})}{e^{-\pi/4}\psi(-e^{-3\pi})}.\tag5.34$$
Putting $\a=\pi/3$ and $\b=3\pi$ in (3.1), we deduce that
$$\frac{\psi(-e^{-\pi/3})}{e^{-\pi/4}\psi(-e^{-3\pi})}=\sqrt3.\tag5.35$$
Using (5.35), we find that
$$PQ=\sqrt3.\tag5.36$$
Using (5.36) in (5.28), we obtain (5.11).
\enddemo
As the pattern of proof of the following theorem is identical with the proof of Theorem 5.4, we skip the proof.

\proclaim{Theorem 5.5} We have
\roster
\item"(i)" $$\frac{\varphi(e^{-\pi/\sqrt5})}{\varphi(e^{-\sqrt5\pi})}=5^{1/4},$$
\item"(ii)" $$\frac{\varphi(e^{-3\pi/\sqrt5})}{\varphi(e^{-3\sqrt5\pi})}
=\sqrt{\frac{10+5\sqrt3+4\sqrt5+2\sqrt{15}}{8+5\sqrt3+4\sqrt5+2\sqrt{15}}},$$
\item"(iii)" $$\frac{\varphi(e^{-\pi/3\sqrt5})}{\varphi(e^{-3\pi/\sqrt5})}
=\frac{9-3\sqrt3+3\sqrt5-\sqrt{15}}{5-3\sqrt3+3\sqrt5-\sqrt{15}},$$
\item"(iv)" $$\frac{\varphi(e^{-\pi/3\sqrt3})}{\varphi(e^{-\pi/\sqrt3})}
=\root4\of{\frac{27+(\root3\of4-1)^4}{3+(\root3\of4-1)^4}},$$
\item"(v)" $$\frac{\varphi(e^{-\pi})}{\varphi(e^{-3\pi})}=\root4\of{6\sqrt3-9},$$
\item"(vi)" $$\frac{\varphi(e^{-\pi/3\sqrt5})}{\varphi(e^{-\sqrt5\pi/3})}
=\sqrt{\frac{10-5\sqrt3+4\sqrt5-2\sqrt{15}}{8-5\sqrt3+4\sqrt5-2\sqrt{15}}},$$
\item"(vii)" $$\frac{\varphi(e^{-\sqrt5\pi/3})}{\varphi(e^{-3\sqrt5\pi})}=
=\frac{9+3\sqrt3+3\sqrt5+\sqrt{15}}{5+3\sqrt3+3\sqrt5+\sqrt{15}},$$
\item"(viii)" $$\frac{\varphi(e^{-\pi/\sqrt3})}{\varphi(e^{-\sqrt3\pi})}=3^{1/4},$$
\item"(ix)" $$\frac{\varphi(e^{-\pi/\sqrt3})}{\varphi(e^{-3\sqrt3\pi})}
=\root4\of{\frac{27+9(\root3\of4-1)^4}{27+(\root3\of4-1)^4}},$$
\item"(x)" $$\frac{\varphi(e^{-\pi/3})}{\varphi(e^{-3\pi})}=\sqrt3,$$
\item"(xi)" $$\frac{\varphi(e^{-\pi/3})}{\varphi(e^{-\pi})}=\root4\of{3+2\sqrt3},$$
\item"(xii)" $$\frac{\varphi(e^{-\pi/\sqrt3})}{\varphi(e^{-3\sqrt3\pi})}=\frac{3\root3\of4}{2+\root3\of4},$$
\item"(xiii)" $$\frac{\varphi(e^{-\pi/3\sqrt3})}{\varphi(e^{-\sqrt3\pi})}=4^{-1/3}(2+\root3\of4).$$
\endroster
\endproclaim
Using (2.7) and (2.8) in Theorem 5.4, we obtain the following values of the cubic continued fraction of Ramanujan.
\proclaim{Theorem 5.6} We have
\roster
\item"" $$V(-e^{-\pi\sqrt5})=\frac{(3-\sqrt5)(\sqrt3-\sqrt5)}{4},$$
\item"" $$V(-e^{-\pi/\sqrt5})=\frac{(\sqrt3+\sqrt5)(\sqrt5-3)}{4},$$
\item"" $$V(-e^{-\pi/\sqrt3})=\frac{-1}{\root3\of4},$$
\item"" $$V(-e^{-\pi\sqrt3})=\frac{1-\root3\of4}{2+\root3\of4},$$
\item"" $$V(-e^{-\pi/3\sqrt3})=\frac{-1}{\root3\of{3^{-1}(\root3\of4-1)^4+1}},$$
\item"" $$V(-e^{-\pi})=\frac{1-\sqrt3}{2},$$
\item"" $$V(-e^{-\pi/3})=\frac{-1}{\root3\of{2(\sqrt3-1)}}.$$
\endroster
\endproclaim
\noindent
{\bf Remark.}  A different proof of Theorem 5.6 (i) can be found in [12] and Theorem 5.6 (iii) was first proved
by C. Adiga et al. [2].  Other values of $V$ in Theorem 5.6 appear to be new to literature.

\vskip3mm
\noindent
{\bf Open Problem.} Using (2.10) on can express $V(q^3)$ in terms of $V(q)$ which gives a triplication formula for $V(q).$
Triplication formula is important in the development of elliptic functions to alternative bases.
The only triplication formula known so far is that of the Borweins.  Can one develop a cubic theory associated with $V(q)$
similar to that of Borweins?

\head\bf Acknowledgement \endhead
Authors are thankful to the referee for  his useful suggestions and comments.  This work was supported by
Korea Research Foundation Grant (KRF-2002-050-C00001).

\enddocument